%
%
%
\documentstyle{amsppt}

\rightheadtext{Infinitely differentiable functions}

\define\C{{\Cal C}}

\define\reals{{\Bbb R}}

\define\nats{{\Bbb N}}

\define\dual{{^{\textstyle*}}}
\define\norm#1{{\|#1\|}}
\define\inorm#1{{\|#1\|_{\infty}}}
\define\half{{\textstyle\frac12}}
\define\third{{\textstyle\frac13}}

\define\seventh{{\textstyle\frac17}}

\define\ball{\text{\rm ball}\,}

\define\less{\setminus}

\define\<{\langle}
\define\>{\rangle}

\topmatter
\title Smooth functions and partitions of unity on certain Banach spaces 
\endtitle
\author Richard Haydon \endauthor
\address{Brasenose College, Oxford, OX1 4AJ, UK}\endaddress
\email{richard.haydon\@brasenose.oxford.ac.uk}\endemail
\endtopmatter

\document

\define\DGZ{2}
\define\Hajek{3}
\define\DevFonHaj{1}
\define\Haydon{4}
\define\Haydontwo{5}
\define\Torunczyk{6}

\heading
Introduction
\endheading

In an earlier paper [\Haydon], the author sketched a method, based on the use 
of ``Talagrand operators'', for defining infinitely differentiable equivalent 
norms on the spaces $\C_0(L)$ for certain locally compact, scattered spaces 
$L$. A special case of this result was that a $\C^\infty$ renorming exists on 
$\C_0(L)$ for every countable locally compact $L$. Recently, H\'ajek [\Hajek] 
extended this result by showing that a real normed space $X$ admits a 
$\C^\infty$ renorming whenever there is a countable subset of the unit ball of 
$X\dual$ on which every element of $X$ attains its norm, that is to say, a 
countable {\sl boundary}. This suggested to the author that the locally 
compact topology on $L$ was perhaps not essential in [\Haydon], and in the 
first part of the present paper we shall develop the methods of that work in 
a way that does not require such a topology. We obtain infinitely 
differentiable norms on certain (typically non-separable) Banach spaces $X$ as 
well as on some certain injective tensor products $X\otimes_\epsilon E$. 

In the second part of the paper we present a lemma about partitions of unity. 
It is an open problem whether a non-separable Banach space with a $\C^k$ norm 
(or, more generally, a $\C^k$ ``bump function'') admits $\C^k$ partitions of 
unity, though many partial results in this direction are known. Our lemma 
enables us to show that the answer is yes for classes of Banach spaces that 
admit projectional resolutions of the identity. In particular, we show that 
the space $\C_0[0,\Omega)$ admits $\C^\infty$ partitions of unity for every  
ordinal $\Omega$. Results from this paper are used in [\Haydontwo] to give 
examples of Banach spaces admitting infinitely differentiable bump functions 
and partitions of unity but no smooth norms.   

Our notation and terminology are mostly standard and, whenever possible, we 
have followed the conventions of [\DGZ]. Although that work contains 
everything that the reader will need in order to understand the present paper, 
we recall for convenience a few facts and definitions. It should be noted that 
we are concerned only with real, as opposed to complex, Banach spaces. When we 
refer to a function on a Banach space as being of class $\C^k$, where $k$ is a 
positive integer, it is the standard (Fr\'echet) notion of smoothness that we 
are employing. Making a mild abuse of language, we shall say that a  norm 
$\norm\cdot$ on a Banach space $X$ is of class $\C^k$ if the function 
$x\mapsto \norm x$ is of that class {\it on the set }$X\less\{0\}$. (Of 
course, no norm is differentiable at $0$.) 

A {\sl bump function} on a Banach space $X$ is a function $\phi :X\to \reals$ 
with bounded, non-empty support. On finite-dimensional spaces, $\C^\infty$ 
bump functions are plentiful (and fundamental to the theory of distributions). 
The existence of a $\C^1$ bump function on an infinite-dimensional Banach 
space $X$ is already a strong condition. For a separable space $X$, it is 
equivalent to separability of the dual space $X\dual$, and to the existence of 
an equivalent $\C^1$ norm. More generally, the existence of a $\C^1$ bump on 
$X$ implies that $X$ belongs to the important class of {\it Asplund spaces.} 
Whether every Asplund space admits a $\C^1$ bump is an open problem, as is the 
relationship between existence of a $\C^k$ bump and of a $\C^k$ norm on a 
separable space once $k$ is greater than 1.  On the other hand, the 
equivalences in the following proposition are very easy to establish. 

\proclaim{Proposition 1}
For a real Banach space $X$ and $k\in \nats\cup\{\infty\}$, the following are 
equivalent: 
\roster
\item $X$ admits a $\C^k$-bump function;
\item there exists a real number $R>1$ and a function $\psi:X\to \reals$, of 
class $\C^k$, such that $0\le \psi\le 1$, $\psi(x)=0$ when $\norm x\le 0$ and 
$\norm x = 1$ when $\norm x\ge R$; 
\item there is a function $\theta:X\to \reals$, of class $\C^k$, such that 
$\theta(x)\to \infty$ as $\norm x\to \infty$.
\endroster
\endproclaim
\demo{Proof}
(1)$\implies$(2):

Let $\phi$ be a $\C^k$ bump function with $\phi(0)=1$. There exist positive real numbers $\delta$ and $M$ such that $\phi(x)\ge \frac23$ when $\norm x\le \delta$ and $\phi(x)=0$ when $\norm x\ge M$. Let $\pi:\reals\to [0,1]$ be a $\C^\infty$ function with $\pi(t)=0$ for $t\ge \frac23$ and $\pi(t)=1$ for $t\le \third$. Then the function $\psi$ defined by 
$$
\psi(x) = \pi(\phi(\delta x))
$$
has the required properties, with $R=\delta^{-1}M$.

(2)$\implies$(3):

Given $R$ and $\psi$ as in (2), we may define
$$
\theta(x) = \sum _{n=0}^\infty \psi(R^{-n}x),
$$
noting that on each ball $\{x\in X: \norm x < N\}$ the sum has only finitely many nonzero terms.  

(3)$\implies$(1):

Given $\theta$, we define $\phi(x)=\pi(\theta(x)-\theta(0))$, where $\pi:\reals\to\reals$ is the function already used above.
\qed\enddemo

Many of our results concern spaces that are subspaces of the space 
$\ell_\infty(L)$ of bounded real-valued functions on a set $L$. For elements 
$f$ of $\ell_\infty(L)$ we use ``coordinate'' notation, writing $(f_t)_{t\in 
L}$ and thinking of the $f_t$ as coordinates. A certain class of very
nice functions, already well-established in the literature, will be of 
particular importance. We shall say that a function $\phi $, defined on a 
subset $D$ of $\ell_\infty(L)$, {\it depends locally on finitely many
coordinates} if, for each $f^0$ in $D$, there exist an open neighbourhood $G$ 
of $f^0$ in $D$ and a finite subset $M$ of $L$ such that, for $f\in G$, the 
value of $\phi (f)$ depends only on $f_t$ $(t\in M)$.

\heading
Infinitely differentiable norms
\endheading

\proclaim{Theorem 1}
Let $L$ be a set and let $U(L)$ be the subset of the direct sum $\ell_\infty(L)
\oplus c_0(L)$ consisting of all pairs $(f,x)$ such that $\inorm{f}$ and 
$\inorm{x}$ are both strictly less than $\inorm{|f|+\half |x|}$. The space 
$\ell_\infty(L) \oplus c_0(L)$ admits an equivalent norm $\norm{\cdot}$ with 
the following properties: 
\roster
\item $\norm{\cdot}$ is a lattice norm, in the sense that $\norm{(g,y)}\le 
\norm{(f,x)}$ whenever $|g|\le |f|$ and $|y|\le |x|$;
\item $\norm{\cdot}$ is infinitely differentiable on the open set $U(L)$;
\item locally on $U(L)$, $\norm{(f,x)}$ depends on only finitely many non-zero 
coordinates; that is to say, for each $(f^0,x^0)\in U(L)$ there is 
a finite $N\subseteq L$ and an open neighbourhood $V$ of $(f^0,x^0)$ in $U(L)$, 
such that for $(f,x)\in V$ the norm $\norm{(f,x)}$ is determined by the values 
of $f_t$ and $x_t$ with $t\in N$ and such that $f_t\ne0$, $x_t\ne 0$ for all 
such $(f,x)$ and $t$. 
\endroster \endproclaim 

We start the proof of Theorem~1 by fixing a strictly increasing 
$\C^\infty$ function $\varpi :[0,2)\to [0,\infty)$ such that 
$\varpi(u)\to\infty$ as $u\uparrow 2$ and $\varpi(u)=0$ for $u\le 1$. The 
inverse function $\varpi^{-1}$ is $\C^\infty$ and strictly increasing from 
$(0,\infty)$ onto $(1,2)$. We define $\theta :[0,\infty\to[0,\infty)$ by 
$$ 
\theta(c)=\int_0^c\frac{\text dv}{\varpi^{-1}(v)},
$$ 
and start by recording some facts about this function. The easy proofs are 
left to the reader.

\proclaim{Lemma 1}
\roster
\item The function $\theta $ is strictly increasing and strictly concave  from 
$[0,\infty)$ onto $[0,\infty)$. It is of class $\C^\infty$ on $(0,\infty)$, 
with $\theta '(c) = 1/\varpi^{-1}(c)$ $(c>0)$, and differentiable at 0 with  
$\theta'(0)=\lim_{c\downarrow 0}1/\varpi^{-1}(c)=1$. 
\item The composite function $\theta \circ\varpi:[0,2)\to [0,\infty)$ is 
infinitely differentiable, with 
$$
(\theta \circ\varpi)'(u)=\left\{\matrix u^{-1}\varpi'(u) &\qquad (u\ge 0)\\
                                0 & \qquad        (u=0)
                        \endmatrix\right.
$$

\item We have $\half c<c\theta'(c)<\theta (c)<c$ for all positive $c$. 
\endroster\endproclaim 

The next lemma can be regarded as the finite-dimensional part  of the proof of 
Theorem~1. 

\proclaim{Lemma 2}
Let $N$ be a finite set, let $\eta $ be a positive real number and and let $W$ 
be the subset of $\reals^N\times \reals^N$ 
consisting of all $(\bold f,\bold x)$ such that $\inorm{|\bold f|+\half 
|\bold x|}>\max\{\inorm{f},\inorm{x}\}+\eta $. Let the functions 
$F:\reals^N\times\reals^N\times[0,\infty)^N\to \reals$,
$G:\reals^N\times\reals^N\to \reals$ be defined by
$$
\align
F(\bold f,\bold x,\bold c)&= \exp\biggl[-\sum_{t\in N} c_t\biggr] \sum_{t\in N} 
\biggl[ c_t|f_t|+\theta (c_t)|x_t|\biggr]\\
G(\bold f,\bold x) &= \sup_{\bold c\in [0,\infty)^N} F(\bold f,\bold x, \bold 
c).
\endalign
$$
If $(\bold f,\bold x)\in W$ then the supremum in the definition of $G(\bold 
f,\bold x)$ is attained at a unique $\bold c$; this $\bold c$ has the property 
that $c_t=0$ whenever $|f_t|\le \eta $ or $|x_t|\le \eta $. The 
function $G$ is of class $\C^\infty$ on $W$. 
\endproclaim
\demo{Proof}
To start with, let us consider a fixed $(\bold f,\bold x)\in W$. We have
$$
\align
\frac{\partial F}{\partial c_s}&= \exp\biggl[-\sum_{t\in N} 
c_t\biggr]\Bigl[|f_s|+\theta'(c_s)|x_s| -\sum_{t\in N} \bigl[c_t|f_t|+\theta 
(c_t)|x_t|\bigr] \Bigr]\\         
                               &\le\exp\biggl[-\sum_{t\in N} 
c_t\biggr]\Bigl[(1-c_s)|f_s|+(\theta'(c_s)-\theta (c_s))|x_s|\Bigr],
\endalign
$$
which is non-positive when $c_s\ge 1$. Thus the supremum in the definition of 
$G(\bold f,\bold x)$ may be taken over the compact set $[0,1]^N$ rather than 
over $[0,\infty)^N$; this supremum is thus attained at some $\bold c$. By 
elementary calculus, any $\bold c$ at which the supremum is attained 
satisfies, for all $s$, either
$$
\align
c_s&>0 \quad\text{ and }\quad|f_s|+\theta '(c_s)|x_s|= \nu \tag 1a \\
\text{ or }\quad c_s&=0\quad \text{ and }\quad|f_s|+\theta'(0)|x_s|\le\nu\tag 1b
\endalign
$$
where $\nu =    \sum_{t\in N} \biggl[c_t|f_t|+\theta (c_t)|x_t|\biggr] $. 
In case (1a) we have 
$$
 c_s = \varpi\Bigl(\frac{|x_s|}{\nu -|f_s|}\Bigr) \tag2
$$
because $\varpi$ is the function inverse to $1/\theta '$. In fact, this 
equality holds in case 1b as well because then $|f_t|+|x_t| = 
|f_t|+\theta'(0)|x_t| \le \nu $, whence 
$|x_t|/(\nu -|f_t|)\le 1$ and $\varpi(|x_t|/(\nu -|f_t|))=0=c_t$. 
Thus $\nu $ is a solution of 
$$ 
\nu = \sum_{t\in N} 
\biggl[\varpi\Bigl(\frac{|x_t|}{\nu -|f_t|}\Bigr)|f_t|+ 
\theta\circ\varpi\Bigl(\frac{|x_t|}{\nu -|f_t|}\Bigr)|x_t|\biggr] \tag3
$$    
Since the right hand side of this equation is a decreasing function of $\nu $ 
it has only one solution. By equation (2), we now see that $c_s$ $(s\in N)$ are 
uniquely determined too. 

Because $(\bold f,\bold x)\in W$, there is some $s$ such that 
$|f_s|+\half|x_s|>\max\{\inorm{f},\inorm{f}\}+\eta $; since $\theta '
(c_s)>\half $ we have $\nu \ge|f_s|+\half|x_s|$ by (1a) and (1b). Thus $\nu 
>|f_t|+\eta $ and $\nu >|x_t|+\eta $ for all $t$. So if either $|x_t|\le \eta 
$ or $|f_t|\le \eta $ it must be that (1b) holds, with $c_t=0$.

We now move on to consider the behaviour of $\nu=\nu(\bold f,\bold x) $ and of 
$c_t=c_t(\bold f,\bold x)$ as $(\bold f,\bold x)$ varies over $W$. We consider 
the function $H$ defined on the open set $V=\{(\bold f,\bold x,\nu ):(\bold 
f,\bold x)\in W\text{ and }\nu >\max\{\inorm f,\inorm x\}+\eta \}$ of 
$\reals^N\times \reals^N\times\reals$ by 
$$
H(\bold f,\bold x,\nu )=\nu -\sum_{t\in N} 
\Bigl[\varpi\Bigl(\frac{|x_t|}{\nu -|f_t|}\Bigr)|f_t|+ 
\theta\circ\varpi\Bigl(\frac{|x_t|}{\nu -|f_t|}\Bigr)|x_t|\Bigr].
$$    
We have already noted that for each $(\bold f,\bold x)\in W$ there is a unique 
$\nu=\nu (\bold f,\bold x) $ such that $H(\bold f,\bold x,\nu )=0$. It is also 
easy to verify that $\frac{\partial H}{\partial \nu }\ge 1$ everywhere on $V$. 
Thus the infinite differentiability of $(\bold f,\bold x)\mapsto \nu (\bold 
f,\bold x)$ will follow from the Implicit Function Theorem provided we can 
show that $H$ is itself infinitely differentiable. The absolute value signs 
appear to present a problem on a neighbourhood of a point where one of the 
variables $f_t$ or $x_t$ is zero. However, as soon as either $|f_t|$ or 
$|x_t|$ is smaller than $\eta $, the terms $\varpi\Bigl(\frac{|x_t|}{\nu-
|f_t|}\Bigr) |f_t| $ and $\theta\circ\varpi\Bigl(\frac{|x_t|}{\nu -
|f_t|}\Bigr)|x_t|$ vanish, showing that we do not have a problem at all. 

Once we have shown that $\nu $ varies in an infinitely differentiable fashion 
with $(\bold f,\bold x)$, it follows from (2) that the same is true for all 
the $c_t$ and hence for $G$. \qed
\enddemo

We now take up the proof of the theorem. 
We define a norm $\norm{\cdot}$ on $\ell_\infty(L)\oplus c_0(L)$ by
$$
\align
\norm{(f,x)} = \sup\Bigl\{\text e^{-\sum_{t\in L} d_t}\sum_{t\in L}&\biggl[ 
d_t|f_t|+\theta (d_t)|x_t|\biggr]:\ d_t\ge 0\text{ for all }t \text{ and }\\
&\ d_t =0\text{ for all but finitely many }t \Bigr\}. 
\endalign
$$
and claim that this has the properties we require. It is clear that 
$\norm{\cdot}$ is a lattice norm and that 
$$
\text e^{-1}\max\left\{\inorm{f},\half\inorm{x}\right\}\le \norm{(f,x)}\le 
\text{e}^{-1}\left(\inorm{f}+\inorm{x}\right).
$$

For $(f,x)\in U(L)$ we define 
$$
\align
\xi (f,x) &= \inorm{|f|+\half|x|}-\max\{\inorm{f},\inorm{x}\}\\
M(f,x) &=\{t\in L:|f(t)|+|x(t)|\ge \inorm{|f|+\half|x|}\}\\
N(f,x) &=\{t\in L:|f(t)|+|x(t)|\ge \inorm{|f|+\half|x|}-\half\xi(f,x)\}
\endalign
$$
and note that $N(f,x)$ is a finite set, because $x\in c_0(L)$ and 
$N(f,x)\subseteq \{t:|x_t|\ge\half \xi(f,x) \}$. We shall show first that in 
the definition of $\norm{(f,x)}$ it is enough to take the supremum over 
families $d=(d_t)$ such that $d_t=0$ for all $t\notin M(f,x)$. Indeed, let 
$(f,x)$ be in $U(L)$, and suppose that $d=(d_t)_{t\in L}$ is such that 
$d_{t_1}>0$ for some $t_1\notin M(f,x)$. Let $t_0$ be chosen so that 
$|f(t_0)|+\half|x(t_0)|=\inorm{|f|+\half|x|}$ and let $d'=(d'_t)$ be defined 
by 
$$
d'_t=\left\{\matrix\format \c & \qquad\qquad\l \\ 
        d_t & \text{ if } t\notin \{t_0,t_1\}\\
            0   & \text{ if } t = t_1\\
            d_{t_0}+d_{t_1} &   \text{ if } t=t_0.
     \endmatrix\right.
$$
We note that $\sum_t d'_t=\sum_t d_t$ and that $\theta (d'_{t_0})-
\theta (d_{t_0}) >\half d_{t_1}$, because $\theta '$ is everywhere 
greater than $\half$.  
$$
\align
\sum_{t\in L} \bigl[d'_t|f_t|+\theta (d'_t)|x_t|\bigr]-
&\sum_{t\in L} \bigl[d_t|f_t|+\theta(d_t)|x_t|\bigr]\\   &= 
d_{t_1}\bigl[|f_{t_0}|-|f_{t_1}|\bigr] + (\theta (d'_{t_0})-\theta (d_{t_0}))
|x_{t_0}| - \theta (d_{t_1})|x_{t_1}|\\
&\ge d_{t_1}\bigl[ |f_{t_0}|+\half|x_{t_0}|-|f_{t_1}|-|x_{t_1}|\bigr]
\endalign
$$
and this is strictly positive by our assumptions about $t_0$ and $t_1$. In 
this way, we may reduce to 0 all coordinates $d_t$ with $t\notin M(f,x)$ while 
increasing the value of $\exp\bigl[-\sum_t d_t\bigr]
\bigl[\sum_t d_t|f_{t}|+\theta(t)|x_t|\bigr]$. 
                           
We now set about finding the neighbourhoods $V$ and finite sets $N$ referred 
to in (3). Given $(f^0,x^0)\in U(L)$, we set $N=N(f^0,x^0)$ and define $V$ to be 
the open set 
$$
V= \{(f,x):\ \inorm{f-f^0},\ \inorm{x-x^0} < \seventh\xi (f^0,x^0)\}.
$$
It is easy to check that if $(f,x)\in V$ then $\xi (f,x)>\half\xi(f^0,x^0)$ 
and $M(f,x)\subseteq N$. By what we have already proved, this shows that on 
the open set $V$ our norm depends only on coordinates in the finite set $N$.

Moreover, for $(f,x)\in V$ we have
$$
\norm{(f,x)}= 
\sup_{\bold c\in [0,\infty)^N} F(\bold c,(f_t)_{t\in N},(x_t)_{t\in N})
$$
where $F:[0,\infty)^N\times \reals^N\times \reals^N$ is the function
$$
F(\bold d,\bold f,\bold x)=\text {exp\,}(-\sum_{t\in N}d_t)\sum_{t\in N}
\bigl[ d_t|f_t|+ \theta (d_t)|x_t|\bigr].
$$
We can thus apply Lemma~2 (with $\eta =\half\xi (f^0,x^0))$ to conclude that 
$\norm{\cdot}$ is infinitely differentiable on $V$. \qed

\proclaim{Remark}{\rm
The property of depending locally on finitely many non-zero coordinates is 
useful in applications. For instance, if $G$ is an open subset of 
a Banach space $X$ and $\phi :G\to \ell_\infty(L)\oplus c_0(L)$ is a continuous
mapping, taking values in $U(L)$, then $x\mapsto \norm{\phi (x)}$ is of class 
$\C^k$ provided only that each coordinate map $x\mapsto \phi (x)_t$ $(t\in L)$ 
is of that class {\it on the set where it is non-zero}. This is particularly 
handy when the coordinate maps are themselves norms or seminorms on $X$.}
\endproclaim

In the following corollary to Theorem 1 we use the above remark to deduce a 
renorming result about injective tensor products. We recall some facts about 
such products. If $X$ and 
$E$ are Banach spaces and $\xi ,\eta $ are elements of the dual spaces 
$X\dual, E\dual$ respectively, then a linear form $\xi \otimes \eta $ may 
be defined on the algebraic tensor product $X\odot E$ by
$$
(\xi \otimes \eta )(\sum_{j=1}^n x_j\otimes e_j) = \sum_j \<\xi ,x_j\>
\<\eta ,e_j\>.
$$
The injective tensor product $X\otimes _\epsilon E$ is defined to be the 
completion of $X\odot E$ for the norm defined by 
$$
\norm{z}_\epsilon = \max \{|\<(\xi \otimes \eta ),z\>|: \xi \in \ball X\dual,\ 
\eta \in\ball E\dual,\}.
$$
The elementary tensor forms $\xi \otimes \eta $ extend by continuity to  
$X\otimes _\epsilon E$ and $\{\xi \otimes \eta :\xi \in \ball X,\ \eta \in 
\ball E\}$ is a weak* compact subset of $\ball(X\otimes _\epsilon  E)\dual$ on 
which every element of  $X\otimes _\epsilon E$ attains its norm. 

If $Q:X_1\to X_2$ and $R:E_1\to E_2$ are bounded linear operators then a 
bounded linear operator $Q\otimes R:X_1\otimes_\epsilon  E_1 \to
X_2 \otimes_\epsilon  E_2$ is determined by $(Q\otimes R)(x\otimes e)= 
(Qx)\otimes (Re)$. A special case is the so-called ``slice map'' $I_X\otimes 
\eta :X\otimes E\to X$ derived from an element $\eta $ of $E\dual$ and given 
by $(I_X\otimes \eta) (x\otimes e)= \<\eta ,e\>x$. By our earlier remark about 
the attainment of norms on elementary tensor forms, we see that for any $z\in
X\otimes_\epsilon  E$ there is some $\eta \in \ball E\dual$ with 
$\norm{z}_\epsilon = \norm{(I_X\otimes \eta )(z)}.$

In the special case where $X$ is a space $\ell_\infty(L)$, then 
$X\otimes_\epsilon  E$ identifies isometrically with a subspace of the 
vector-valued function space $\ell_\infty(L;E)$ (the elementary tensor
$(x_t)_{t\in L}\otimes e$ corresponding to $t\mapsto x_t e$). The effect of a 
slice map on $z\in\ell_\infty(L)\otimes_\epsilon E$ regarded as an element
$(z_t)_{t\in L}$ of $\ell_\infty(L;E)$ is simply
$$
(I_{\ell_\infty(L)}\otimes \eta )(z) = (\<\eta ,z_t\>)_{t\in L}.
$$

\proclaim{Corollary 1}
Let $X$ be a Banach space and let $L$ be a set. Suppose that there exist a 
linear homeomorphic embedding $S:X\to\ell_\infty(L)$ and a linear operator 
$T:X\to c_0(L)$ with the property that $(S(x),T(x))$ is in the set $U(L)$ of 
Theorem 1 whenever $x$ is a non-zero element of $X$. Then $X$ admits an 
equivalent $\C^\infty$ norm. Moreover, whenever the Banach space $E$ admits an 
equivalent $\C^k$ norm so does the injective tensor product $X\otimes_\epsilon 
E$. 
\endproclaim 
\demo{Proof} 
It is clear that, using the norm on $\ell_\infty(L)\oplus c_0(L)$ that we 
defined in Theorem~1, we may set
$$
\norm{x} = \norm{(S(x),T(x))}
$$
and obtain an infinitely differentiable norm on $X$. 

Now let $E$ be a Banach space with a $\C^k$ norm $\norm{\cdot}_E$. For $f\in 
\ell_\infty(L;E)$  we define $Nf\in \ell_\infty(L)$ by
$$
(Nf)_t= \norm{f_t}_E.
$$
The operators $S$ and $T$ induce $S\otimes I_E$ and $T\otimes I_E$, taking the 
injective tensor product $X\otimes _\epsilon E$ into $\ell_\infty(L)\otimes 
_\epsilon E$ and $c_0(L)\otimes _\epsilon E$ respectively. Identifying 
$c_0(L)\otimes _\epsilon E$ with $c_0(L;E)$ and $\ell_\infty(L)\otimes 
_\epsilon E$ with a subspace of $\ell_\infty(L;E)$, we may define a norm on 
$X\otimes_\epsilon  E$ by
$$
\norm{z} = \norm{(N((S\otimes I_E)(z)),(N((T\otimes I_E)(z)))}.
$$
The coordinate maps  $x\mapsto(N((S\otimes I_E)(x)))_t = \norm{(S\otimes 
I_E)(x)_t}$ and $x\mapsto (N((T\otimes I_E)(x)))_t = 
\norm{(T\otimes I_E)(x)_t}$ are of class $\C^k$ except where they vanish. 
Thus, by the above remark, we shall be able to conclude that we have a $\C^k$ 
norm on $X\otimes_\epsilon E$ provided we can show that 
$(N((S\otimes I_E)(z)),(N((T\otimes I_E)(z)))$  is in $U(L)$ whenever $z\in 
(X\otimes _\epsilon E)\less\{0\}$. This is not really difficult, being just a 
matter of disentangling tensor notation.

For such a $z$ we choose some $\eta \in\ball E\dual$ with 
$\norm {I_{\ell_\infty(L)}\otimes \eta )((S\otimes I_E)(z))}_\infty = 
\norm{(S\otimes I_E)(z)}_\infty $. We then note that 
$(I_{\ell_\infty(L)}\otimes \eta )\circ(S\otimes I_E) = S\circ(I_X\otimes 
\eta)$. Our hypothesis about $S$ and $T$, applied to $x= (I_X\otimes \eta )(z)$ 
tells us that there is some $t\in L$ with 
$$
|S((I_X\otimes \eta) (z))_t|+\half |T((I_X\otimes \eta) (z))_t| > 
\inorm{S((I_X\otimes \eta) (z))}.
$$
Thus
$$
\align
N((S\otimes I_E)(z))_t&+\half N((T\otimes I_E)(z))_t \  = \  
\norm{(S\otimes I_E)(z)_t}_E+\half \norm{(T\otimes I_E)(z)_t}_E\\
&\ge |\<\eta ,(S\otimes I_E)(z)_t\>|+\half| \<\eta,(T\otimes I_E)(z)_t\>|\\
&= |(I_{\ell_\infty(L)}\otimes\eta)(S\otimes I_E)(z)_t|+\half| 
(I_{c_0(L)}\otimes\eta)(T\otimes I_E)(z)_t|\\
&=|S((I_X\otimes \eta) (z))_t|+\half |T((I_X\otimes \eta) (z))_t| \\
&>\inorm{S((I_X\otimes \eta) (z))}\\
&=\norm{(S\otimes I_E)(z)}_\infty
\endalign
$$
A similar argument involving an $\eta $ chosen so that 
$\norm{(I_{c_0(L)}\otimes \eta )((T\otimes I_E)(z)}$ is equal to
$ \norm{(T\otimes I_E)(z)}$ enables us to finish the proof. \qed \enddemo 

As we may now note, the above corollary allows us to reprove H\'ajek's theorem 
from [\Hajek], though not, of course, the more recent, and very strong, result 
of [\DevFonHaj], according to which any norm on a Banach space with countable 
boundary may be approximated by {\sl analytic} norms. 

\proclaim{Corollary 2 [H\'ajek]} Let $X$ be a Banach space which admits a 
countable boundary. Then $X$ admits a $\C^\infty$ renorming and $X\otimes 
_\epsilon E$ admits a $\C^k$ renorming whenever $E$ does. \endproclaim 
\demo{Proof} Let  $\{\xi _n:n\in\nats\}$ be a countable boundary for $X$ and 
define $S:X\to \ell_\infty(\nats)$ and $T:X\to c_0(\nats)$ by $(Sx)_n=\<\xi 
_n,x\>$, $(Tx)_n=2^{-n}\<\xi _n,x\>$. It is easy to see that $(S(x),T(x))\in 
U(L)$ when $x$ is  a non-zero element of  $X$, since  for any $x$ there exists 
$n$ with $\<\xi _n,x\>=\norm x$. \qed\enddemo 

Extending slightly the terminology of [\Haydon], we shall say that a bounded 
linear operator $T$ from a subspace $X$ of $\ell_\infty(L)$ into $c_0(L)$ is a 
{\it Talagrand operator for $X$} if for every non-zero $x$ in $x$ there exists 
$t\in L$ with $|x_t|=\inorm{x}$ and $(Tx)_t\ne 0$. It is clear that 
Corollary~1 is applicable to any such space, taking $S$ to be the identity 
operator. In the particular case where $L$ is equipped with a locally compact 
topology and $X$ is the space $\C_0(L)$ of continuous functions, vanishing at 
infinity, we retrieve the results of [\Haydon]. The classic example of a 
Talagrand operator is defined on the space $\C_0([0,\Omega )$, where $\Omega $ 
is an ordinal, by
$$
(Tf)_\alpha =f_\alpha -f_{\alpha +1}.
$$
This has the required property since for any non-zero $f\in \C_0([0,\Omega ))$ 
there is a {\it maximal} $\alpha $ with $|f_\alpha |=\inorm f$.
A non-linear version of a Talagrand operator is used in [\Haydon] to give an 
example of  a space admitting a $\C^\infty$ bump function but no smooth norm. 
The author's earlier $\C^1$ version of this result appears as Theorem VII.6.1 
of [\DGZ]. The relevant application of our theorem is the following.

\proclaim{Corollary 3}
Let $X$ be a Banach  space and let $L$ be a set. Suppose that there exist 
continuous mappings $S:X\to \ell_\infty(L)$, $T:X\to c_0(L)$ with the 
following properties: 
\roster
\item for all $x\in X$ the pair $(Sx,Tx)$ is in $U(L)\cup \{0\}$;
\item the coordinates of $S$ and of $T$ are all $\C^k$ functions on the sets 
where they are non-zero;
\item $\inorm{Sx}\to \infty$ as $\norm x\to \infty$.
\endroster                      
Then $X$ admits a $\C^k$ bump function.
\endproclaim\demo{Proof}
Let $\theta :[0,\infty)\to [0,\infty)$ be a $\C^\infty$ function which 
vanishes on $[0,1]$ and which tends to infinity with its argument. The formula 
$$
\phi(x) = \theta (\norm{(Sx,Tx)}   )
$$
defines a $\C^k$  function on $X$ which tends to infinity with $\norm {x}$.
\qed\enddemo

The author does not know whether the results of Corollary~1 about injective 
tensor products extend to the non-linear set-up of Corollary~3. However, in 
the special case of spaces of continuous functions we have the 
following proposition.

\proclaim{Proposition 2}
Let $L$ be a locally compact space such that there exists a function 
$T:\C_0(L)\to c_0(L)$ satisfying
\roster
\item for all $f\in\C_0(L)$ the pair $(f,Tf)$ is in $U(L)\cup \{0\}$;
\item each coordinate of $T$ is a $\C^k$ function depending locally on 
finitely many coordinates.
\endroster                      
Let $E$ be a Banach space admitting a $\C^k$ bump function. Then 
the space $\C_0(L;E)$ also admits such a function. In particular, if each of 
the spaces $L_1,\dots,L_n$ is homeomorphic to an ordinal then 
$\C_0(L_1\times\dots\times L_n;E)$ admits a $\C^k$ bump function.
\endproclaim
\demo{Proof}
Let $\theta $ be a $\C^k$ function on $E$ such that $\theta (x)\to \infty$ as 
$\norm x\to \infty$. For $f\in \C_0(L;E)$ the composition $\theta\circ f$ is in 
$\C_0(L)$ and the pair $ (\theta \circ f,T(\theta \circ f))$ is in 
$U(L)\cup\{0\}$. Moreover, for $t\in L$, the coordinate maps $f\mapsto (\theta 
\circ f)_t$ and $f\mapsto T(\theta \circ f)_t$ are of class $\C^k$ on 
$\C_0(L;E)$.  Let $\rho :\reals\to \reals$ be a $\C^\infty$ function such that 
$\rho (u)=0$ for $u\le 1$ and $\rho (u)\to \infty$ as $u\to \infty$. It is 
easy to check that the function $\phi :\C_0(L;E)\to \reals$ defined by
$$
\phi (f)= \rho \biggl(\norm{(\theta \circ f,T(\theta \circ f))}\biggr)
$$
is of class $\C^k$ and tends to infinity with the norm of its argument. 

When $L$ is an ordinal $\Omega $ (identified with the set $[0,\Omega )$ of 
ordinals smaller than itself), an operator $T$ of the type considered above 
does exist. Indeed, we may use the Talagrand operator
$(Tf)_\alpha = f_\alpha -f_{\alpha +1}$. Thus $\C_0([0,\Omega);E)$ admits a 
$\C^k$ bump function whenever $E$ does. 
Since $\C_0(L_1\times\cdots\times L_n;E)$ may be identified with $\C_0(L_1; 
\C_0(L_2\times\cdots\times L_n;E))$, an easy induction argument finishes the 
proof.
\qed
\enddemo

\heading
Smooth partitions of unity
\endheading

We say that a subset $H$ of a Banach space $X$ admits {\it partitions of unity 
of class} $\C^k$ if, for every open covering $\Cal V$ of $H$, there is a 
locally finite partition of unity on $H$, subordinate to the covering $\Cal 
V$, and consisting of the restrictions to $H$ of functions that are of class 
$\C^k$ on $X$. Once again, the reader is referred to [\DGZ] for more details, 
for the connection between partitions of unity and approximation by smooth 
functions and for Torunczyk's criterion: $H$ admits $\C^k$ 
partitions of unity if and only if there is a $\sigma $-locally finite base 
for the topology of $X$ consisting of $\C^k$-cozero sets (that is to say, sets 
of the form $\{x\in H: \phi (x)\ne 0\}$ with $\phi \in \C^k(X)$). It is not 
known whether $\C^k$ partitions of unity necessarily exist on every space that 
admits a $\C^k$ bump function, though many partial results are known ([\DGZ, 
VIII.3]). The following theorem has hypotheses that are involved but of fairly 
wide applicability.

\proclaim{Theorem 2}
Let $X$ be a Banach space, let $\Gamma $ be a set and let $k$ be a positive 
integer or $\infty$. Let $T:X\to c_0(\Gamma 
)$ be a function such that each coordinate $x\mapsto T(x)_\gamma $ is of class 
$\C^k$ on the set where it is non-zero. For each finite subset $F$ of $\Gamma $,
let $R_F:X\to X$ be of class $\C^k$ and  assume that the following hold:
\roster
\item
for each $F$, the image $R_F[X]$ admits $\C^k$ partitions of unity;
\item
$X$ admits a $\C^k$ bump function;
\item
for each $x\in X$ and each $\epsilon >0$ there exists $\delta >0$ such that
$\norm{x-R_Fx}<\epsilon $ if we set $F =\{\gamma \in\Gamma 
:|(Tx)(\gamma )|\ge \delta \}$.
\endroster
Then $X$ admits $\C^k$ partitions of unity.
\endproclaim
\demo{Proof}
By Torunczyk's Criterion, it is enough to show that there is a $\sigma $-%
locally finite base  for the topology of $X$, consisting of $\C^k$-cozero-%
sets. By hypothesis, each $R_F[X]$ admits a $\sigma $-locally finite base 
$\Cal V_F$ consisting of $\C^k$-cozero-sets. Since $X$ admits a $\C^k$-bump 
function, there is a neighbourhood base of $0$ in $X$ consisting of $\C^k$-%
cozero-sets, say $U_n$ $(n\in\nats)$. We introduce the covering $\Cal W$ of 
$c_0(\Gamma )$ consisting of $W_\emptyset=c_0(\Gamma )$, together with 
all sets
$$
W_{F,q,r}  = \{y\in c_0(\Gamma ): \min_{\gamma\in F}|y(\gamma )|> 
r\text{ and }\sup_{\gamma \in \Gamma \less F}|y(\gamma )|<q\}
$$
with $F$ a finite non-empty subset of $\Gamma $, and $q,r$ positive rational 
numbers with $q<r$. We note that $\Cal W$ is $\sigma $-locally finite, and 
that its members are $\C^\infty$-cozero-sets. 

In $X$ we consider the family of all sets of the form
$$
T^{-1}[W_{F,q,r}]\cap R_F^{-1}[V]\cap (R_F-I)^{-1}[U_m]
$$
with $m$ a positive integer, $F$ a finite subset of $\Gamma $, $q,r$ positive 
rationals with $q<r$ and $V\in \Cal V_F$. It is easy 
to check that this family is a $\sigma $-locally finite family of $\C^k$ 
cozero sets. We have to show that it forms a base for the topology of $X$. 

Let $x$ be in $X$ and let $\epsilon >0$ be given. We fix $m$ so that
$$
U_m\subseteq \third\epsilon\, \ball X,
$$
and, using (3), choose $\delta >0$ so that
$$
x-R_F(x)\in U_m
$$
when we set $F=\{\gamma \in\Gamma :|(Tx)(\gamma )|\ge \delta \}$. Because 
$Tx\in c_0$ there exist rationals $q,r$ with $0<q<r<\delta $ such that 
$|(Tx)(\gamma )| < q$ whenever $\gamma \in \Gamma \less F$. Thus $x$ is in 
$T^{-1}[W_{F,q,r}]$. Since $\Cal V_F$ is a base for the topology of $R_F[X]$, 
there exists $V\in \Cal V_F$ such that 
$$
R_F(x)\in V\subseteq R_F(x)+\third\epsilon\, \ball X.
$$

It follows that $x$ is in the set
$$
T^{-1}[W_{F,n}]\cap R_F^{-1}[V]\cap (R_F-I)^{-1}[U_m].
$$
If $x'$ is any other member of this set, then we have
$$
\norm{R_F(x)-R_F(x')}\le \epsilon /3 
$$
because $R_F(x')\in V$, while
$$
\norm{R_F(x')-x'} \le \epsilon /3,
$$
because $(R_F(x')-x')\in U_m$. Thus  $\norm{x-x'}<\epsilon $, which is what we 
wanted to prove.
\qed\enddemo

It should be noted that the mappings $T$ and $R_F$ in the theorem are not 
assumed to be linear; a non-linear $T$ is used in [\Haydontwo] to give an 
example where $\C^\infty$ partitions of unity exist on a space with no smooth 
norm.   However, the theorem offers some improvements on existing results even 
when only linear operators are involved. A special case of the corollary that 
follows occurs when the $R_\alpha $ form a ``projectional resolution of the 
identity'' on $X$. It is thus a result that is more general, as well as a bit 
simpler to prove, than the implication (vi)$\implies$(v) in Theorem VII.3.2 of 
[\DGZ]. 

\proclaim{Corollary 4}
Let $X$ be a Banach space admitting a $\C^k$ bump function. Let $\Omega $ be 
an ordinal and let $R_\alpha $ $(\alpha <\Omega )$ be a family of $\C^k$ 
functions from $X$ to $X$ having the property that, for every $x\in X$, the 
function $Rx:[0,\Omega ]\to X$ defined by $(Rx)_\alpha =R_\alpha x$ $(\alpha 
<\Omega )$, $(Rx)_\Omega=x$ is continuous. If for each $\alpha $ the image of 
$R_\alpha $ admits $\C^k$ partitions of unity then so does $X$. \endproclaim 
\demo{Proof} 
Since $X$ admits a $\C^k$ bump function there exists a function $\phi :X\to 
[0,1]$, of class $\C^k$ and such that $\phi (x)=0$ on some neighbourhood of 0 
in $X$, while $\phi (x)=1$ whenever 
$\norm{x}\ge 1$. We set $\Gamma =\Omega \times \nats$ and define $T:X\to 
\ell_\infty(\Gamma )$ by 
$$
(Tx)(\alpha ,n)= 2^{-n} \phi (2^n(R_{\alpha +1}x-R_\alpha x )).
$$
By construction, there is some $\eta >0$ such that $\phi (x)=0$ whenever 
$\norm{x}\le \eta $. Given $x\in X$ and  $\epsilon >0$ we fix $m$ such that 
$2^{-m}<\epsilon $ and then note that, because of the continuity of $\alpha 
\mapsto R_\alpha x$, the quantity $\norm{R_{\alpha +1}x-R_\alpha x} $ can 
exceed $2^{-m}\eta $ only for $\alpha $ in some finite set $H$. We thus have 
$|(Tx)_\gamma |\le \epsilon $ except when $\gamma \in H\times \{0,1,2,\dots, 
m-1\}$. This shows that $T$ takes its values in $c_0(\Gamma )$.

To define the ``reconstruction operators'' $R_F$ we        
set $R_\emptyset = R_0$ and $R_F=R_{\alpha (F)+1}$ where, for a finite 
non-empty subset $F$ of $\Gamma $, $\alpha (F)=\max\{\alpha :\exists  n
\text{ with } (\alpha ,n)\in F\}$.  We shall show that Condition (3) 
of Theorem~2 is satisfied. Given $x\in X$ and $\epsilon >0$, it may be that 
$\norm{x-R_\alpha x}<\epsilon $ for all $\alpha <\Omega $; in this case there 
is clearly no problem.  Otherwise, by the continuity of $\alpha \mapsto 
R_\alpha x$ on $[0,\Omega ]$, there is a maximal $\beta <\Omega $ with 
$\norm{x-R_\beta x}\ge \epsilon $.  Again by the continuity of $\alpha 
\mapsto R_\alpha x$, we know that there is some $\gamma >\beta $ such that 
$\norm{R_{\gamma +1}x-R_\gamma x}$ takes a strictly positive value, $\eta $ 
say.  Now we fix $n$ such that $2^n\eta \ge 1$, noting that $(Tx)(\gamma 
,n)=2^{-n}$, and set $\delta =2^{-n}$.  If $F$ is the set $\{(\alpha,m)\in
\Omega \times\nats:|(Tx)(\alpha ,m)|\ge \delta \}$ then $(\gamma ,n)\in F$ 
and so $\alpha (F)\ge \gamma >\beta $, whence $\norm{x-R_Fx}<\epsilon $, as 
required. \qed\enddemo

\proclaim{Corollary 5}
Let $\Omega $ be an ordinal and let $E$ be a Banach space admitting $\C^k$ 
partitions of unity. Then the space $\C([0,\Omega ];E)$ also admits $\C^k$ 
partitions of unity.
\endproclaim
\demo{Proof}          Proceeding by transfinite induction, we may suppose that 
$\C([0,\gamma ];E)$ admits $\C^k$ partitions of unity whenever $\gamma $ is 
an ordinal smaller than $\Omega $. If we define $R_\gamma :\C([0,\Omega ];E)
\to\C([0,\Omega ];E)                $ by 
$$
(R_\gamma f)_\beta =\left\{\matrix f_\beta & \quad (\beta \le \gamma )\\
                                    f_\gamma &\quad (\beta > \gamma )
                            \endmatrix\right.
$$
then the range of $R_\gamma $ is isomorphic to $\C([0,\Gamma ];E)$ and so 
admits $\C^k$ partitions of unity. Moreover, the continuity hypothesis in the 
preceding corollary is certainly satisfied, so that the proof will be 
finished if we know that $\C([0,\Omega ];E)$ admits a $\C^k$ bump function. 
This is true by Proposition~2, since $\C([0,\Omega ];E)$ is isomorphic to 
$X=\C_0([0,\Omega );E)\oplus E$.
\qed\enddemo

\enddocument